\documentclass[12pt]{amsart}
\usepackage{amsthm,amsmath,amssymb,amsthm,amscd}
\usepackage[T1]{fontenc}
\usepackage[utf8]{inputenc}
\usepackage[english]{babel}
\usepackage{times}
\usepackage{fullpage}
\usepackage{enumerate}
\usepackage{xfrac}
\usepackage{lmodern}

\usepackage{hyperref}
\usepackage{backref}

\usepackage{cleveref}
\usepackage{graphicx}
\usepackage{tikz-cd}

\usetikzlibrary{cd}
\usetikzlibrary{positioning}
\usetikzlibrary{matrix}

\newcommand{\B}{\text{B}}                   
\newcommand{\C}{\mathbb{C}}            
\newcommand{\R}{\mathbb{R}} 
\newcommand{\Q}{\mathbb{Q}}
\newcommand{\Z}{\mathbb{Z}}
\newcommand{\N}{\mathbb{N}}
\newcommand{\CP}{\mathbb{C}\mathbb{P}}

\newtheorem{theorem}{Theorem}[section]

\newtheorem{proposition}{Proposition}[section]

\theoremstyle{remark}

\title[Approximation results for compact Vaisman manifolds]{Approximation results for compact Vaisman manifolds}

\author{Daniele Angella}
\address{Dipartimento di Matematica e Informatica “Ulisse Dini”, Universit\`a degli Studi di Firenze, viale Morgagni 67/A, 50134 Firenze, Italy}
\email{daniele.angella@unifi.it}
\email{daniele.angella@gmail.com}

\author{Marco Miceli}
\address{Dipartimento di Matematica e Informatica “Ulisse Dini”, Universit\`a degli Studi di Firenze, viale Morgagni 67/A, 50134 Firenze, Italy}
\email{marco.miceli@edu.unifi.it}
\email{marcoo.miceli@gmail.com}

\author{Giovanni Placini}
\address{Dipartimento di Matematica e Informatica, Universit\`a degli Studi di Cagliari, via Ospedale 72, 09124 Cagliari, Italy}
\email{giovanni.placini@unica.it}

\date{\today ; {\copyright \ D.~Angella, M. Miceli and G.~Placini 2024}}

\subjclass[2020]{53C55; 32Q40;  53C18}

\keywords{locally conformally K\"ahler metric; Vaisman geometry; Hopf manifold; Tian approximation}

\thanks{
DA was partially supported by project PRIN2022 ”Real and Complex Manifolds: Geometry and Holomorphic Dynamics” (code 2022AP8HZ9) and by GNSAGA of INdAM.
GP was partially supported by INdAM and GNSAGA and by ISI-HOMOS - Funded by Fondazione di Sardegna. GP acknowledges financial support by PNRR e.INS (CUP F53C22000430001, codice MUR ECS00000038)}

\dedicatory{
The authors contribute with solidarity and support to projects that aid colleagues and students in Gaza and in all war areas:\\
\href{https://umi.dm.unibo.it/gruppi-umi-2/gruppo-umi-progetti-di-cooperazione-internazionale-con-il-sud-globale}{https://umi.dm.unibo.it/gruppi-umi-2/\\gruppo-umi-progetti-di-cooperazione-internazionale-con-il-sud-globale}.
}

\begin{document}

\begin{abstract}
We extend the Tian approximation theorem for projective manifolds to a class of complex non-K\"ahler manifolds, the so-called Vaisman manifolds. More precisely, we study the problem of approximating compact regular, respectively quasi-regular, Vaisman metrics by metrics induced by immersions, respectively embeddings, into Hopf manifolds.
\end{abstract}
	
\maketitle

\section*{Introduction}

A classical result by Kodaira states that compact complex manifolds with an ample line bundle $L$ can be embedded into a projective space via sections of powers $L^k$ for $k$ large enough.
Generally the Kodaira embedding is {\em not} an isometry.
However, Tian \cite{tian1990set} proved that suitable rescaling of the metrics induced by Kodaira embeddings $\mathcal{C}^2$-approximate the K\"ahler metric as the power $k$ goes to infinity. 
This was later proved independently by Ruan \cite{ruan1998canonical} and Zelditch \cite{zelditch1998szego} to be a smooth approximation.

Both immersion and approximation results have been studied in other settings, including CR manifolds and other geometric structures, Sasakian manifolds, K\"ahler orbifolds, see {\itshape e.g.} \cite{ma-marinescu-Book, ornea-verbitsky-MathAnn, ross-thomas-2, loi-placini} and references therein. 

In this note, we deal with the problem of approximating Hermitian non-K\"ahler metrics by holomorphic embeddings into model spaces. More precisely, we consider the class of compact Vaisman manifolds, the model space being Hopf manifolds.

A {\em Vaisman} structure on a complex manifold $M$ is given by a Hermitian metric $\omega$ that locally admits a conformal change to a K\"ahler metric, with the further property that the local conformal changes define a closed parallel $1$-form $\theta$.
We refer to \cite{ornea-verbitsky-Book} for an up-to-date account of locally conformally K\"ahler and Vaisman geometry.
Vaisman geometry is strictly related to K\"ahler and Sasakian geometries, see \Cref{sec:preliminaries}. In particular, so-called {\em regular} Vaisman manifolds are given by elliptic fibrations over a projective manifold, respectively {\em quasi-regular} Vaisman manifolds have the same structure over a projective orbifold.

The main example of Vaisman manifold is given by the {\em Hopf manifold}, namely a quotient of $\mathbb C^n\setminus\{0\}$ by the free action generated by a holomorphic contraction.
Hopf manifolds play the role of projective spaces in K\"ahler geometry.
For example, it is known \cite{ornea-verbitsky-MathAnn,ornea2010locally} that Vaisman manifolds are embeddable into Hopf manifolds.

In light of this, our first result is the following:

\begin{theorem}\label{thm:a}
Let $(M,\omega,\theta)$ be a compact regular Vaisman manifold. Then there exists a sequence of holomorphic immersions $\varphi_k \colon M \to \mathcal{H}_{N_k+1}$ into semisimple regular Hopf manifolds endowed with Vaisman metrics such that suitable homotheties of the induced structures converge to $(\omega,\theta)$ in the $\mathcal{C}^\infty$-norm.
\end{theorem}

The proof is based on the Tian approximation theorem,
by exploiting the fact that, when the Lee class $[\theta]\in H^1(M,\mathbb{R})$ is a multiple of a rational class, $M$ is covered by a negative line bundle $L^{-1}$ over a projective manifold $X$, that is, the leaf space of the canonical foliation induced by $\theta$ and $J\theta$.
This is combined with the smooth approximation of a generic Vaisman structure by Vaisman structures with the property that $[\theta]$ is a multiple of a rational class.

In \cite[Theorem 2.11]{ornea-verbitsky-MathAnn}, it is proven that any compact quasi-regular Vaisman manifold admits a holomorphic immersion into a classical regular Hopf manifold.
Notice that quasi-regular Vaisman Hopf surfaces cannot be approximated by metrics induced by classical Hopf manifolds. One such example was presented by the first-named author and Zedda in \cite[Proposition 3.10]{angella2019isometric}.
This, together with the fact that a Vaisman submanifold of a regular Vaisman manifold needs to be regular, shows that, in order to get an approximation result for quasi-regular Vaisman structures, one needs to consider embeddings into quasi-regular Hopf manifolds.
Indeed, our second result on quasi-regular Vaisman structures reads

\begin{theorem}\label{thm:b}
Let $(M,\omega,\theta)$ be a compact quasi-regular Vaisman manifold. Then, for any given $m\in\mathbb N$, there exists a sequence of Vaisman structures induced by holomorphic embeddings $\varphi_k \colon M \hookrightarrow \mathcal{H}_{N_k+1}$ into semisimple quasi-regular Hopf manifolds which converge to $(\omega,\theta)$ up to suitable homotheties in the $\mathcal{C}^m$-norm.
\end{theorem}

As in the proof of \Cref{thm:a}, we first reduce to the case of Vaisman structures whose Lee class $[\theta]\in H^1(M,\mathbb{R})$ is a multiple of a rational class. Such quasi-regular structures are obtained as quotients of negative line orbibundle over projective orbifolds.
The proof is now based on the orbifold embedding theorem for polarized cyclic orbifolds into weighted projective spaces by Ross and Thomas \cite{ross-thomas-2}.
Note that this is indeed different from \cite[Theorem 2.11]{ornea-verbitsky-MathAnn} which relies on a classical result of Baily \cite{baily} realizing polarized orbifolds as projective algebraic subvarieties of the projective space. Our approach allows us to get an embedding and to control the induced metrics by a $C^0$-convergence result of \cite{ross-thomas-1, ross-thomas-2}, later improved to a $C^m$-convergence by Loi and the third-named author \cite{loi-placini}.

Concerning the metric approximation in the irregular case, notice that the characteristic foliation of a compact Vaisman manifold $(M,J,\theta,\omega)$ depends only on the complex structure $J$.
Therefore one cannot approximate irregular Vaisman manifolds by holomorphic embeddings into quasi-regular Hopf manifolds as they would be forced to be quasi-regular.
Furthermore, contrary to the Sasakian case, it does not hold that quasi-regular Vaisman structures are dense in the space of Vaisman structures on the smooth manifold $M$.

\subsection*{Structure of the paper} The paper is organized as follows. In \Cref{sec:preliminaries}, we give a brief presentation on basic definitions and known results in Vaisman geometry. \Cref{sec:regular} and \Cref{sec:quasi-regular} are dedicated to the proofs of \Cref{thm:a} and \Cref{thm:b} respectively.

\subsection*{Acknowledgments}
This note is based on the Master thesis of the second-named author.
The authors would like to warmly thank Nicolina Istrati, Liviu Ornea, Sergio Vessella, Victor Vuletescu, and Michela Zedda.

\section{Preliminaries on Vaisman geometry}
\label{sec:preliminaries}
In this section we recall some know results on LCK and Vaisman manifolds. We refer the reader to \cite{ornea-verbitsky-Book} for a comprehensive and up-to-date exposition of developments in LCK geometry, and to \cite{boyer-galicki-Book} for the classical results on Sasakian manifolds. 
\subsection{Locally conformally K\"ahler structures}

A \emph{locally conformally K\"{a}hler} (LCK) manifold is a Hermitian manifold $(M,J,g)$ such that for any point $p$ there exists an open neighbourhood $U_p$ and a $\mathcal{C}^\infty$ function $\varphi_p \colon U_p \to \R$ such that the local metric $g_p = e^{-\varphi_p} \left.g\right|_{U_p}$ is K\"{a}hlerian. This is equivalent to ask the fundamental form $\omega(\cdot,\cdot):=g(\cdot,J\cdot)$ to satisfy the following cohomological condition \begin{equation*} d\omega = \theta \wedge \omega, \end{equation*} where $\theta$ is a closed $1$-form, called the \emph{Lee form}, which encodes the local conformal changes.

If $\theta = d\phi$ for some $\mathcal{C}^\infty$ function $\phi$ on $M$, then $e^{-\phi}\omega$ is a K\"ahler form on $M$: in this case, the manifold $M$ is called \emph{globally conformally K\"{a}hler} (GCK). Since, for our purposes, we are interested in manifolds not admitting K\"ahler metrics, we will always assume the Lee form $\theta$ to be non-exact.
More precisely, we notice that the LCK property is conformally invariant. Indeed, if $g$ is LCK with Lee form $\theta$, then $e^fg$ is LCK with Lee form $\theta + df$. It follows that we can associate to the conformal class of a LCK structure $(\omega,\theta)$ on $M$ a unique cohomological class $[\theta] \in H^1_{dR}(M,\R)$.

Any LCK manifold admits a K\"ahler covering $\pi \colon \widetilde{M} \to M$ on which the deck transformation group $\Gamma:= \frac{\pi_1(M)}{\ker \rho}$ acts by homotheties, where $\rho$ is a group homomorphism of the fundamental group of $M$ into the additive group $\R$.
Thus, $\Gamma \simeq \Z^k$ where $k$ is the {\em LCK rank} of the LCK structure, namely, the dimension of the smallest linear subspace $W \subset H^1_{dR}(M,\Q)$ such that $W \otimes_\Q \R$ contains the cohomology class $[\theta]$ of the Lee form.
If the K\"ahler cover $(\widetilde{M}, \widetilde{\omega})$ admits a K\"ahler potential $\varphi$ that is \emph{automorphic} with respect to the action of the deck transformation group $\Gamma$, {\itshape i.e.} $\gamma^* \widetilde{\omega} = e^{-\rho(\gamma)} \widetilde{\omega}$ for any $\gamma\in\Gamma$, then the corresponding LCK manifold $(M,\omega, \theta)$ is called a \emph{LCK manifold with potential}. In \cite[Claim 2.4]{ornea-verbitsky-JGP}, it is proven that the potential $\varphi$ is proper if and only if the LCK rank of $M$ is equal to $1$.

The first compact examples of LCK metrics were constructed on Hopf manifolds. A {\em Hopf manifold} is a compact complex non-K\"ahler manifold, diffeomorphic to $S^1 \times S^{2n-1}$ with $n \ge 2$, defined as a quotient of $\mathbb C^n\setminus \{0\}$ by the infinite cyclic group $\mathbb Z$ generated by a holomorphic contraction, namely a smooth map $\gamma \colon \mathbb C^n \to \mathbb C^n$ such that a sufficiently big iteration maps any given compact subset of $\mathbb C^n \setminus \{0\}$ onto an arbitrarily small neighbourhood of $0$. When $\gamma$ is linear, we call the corresponding quotient a {\em linear Hopf manifold}. In the same vein, the terms {\em semisimple} or {\em diagonal} Hopf manifolds will also be used.

\subsection{Vaisman structures}

A LCK manifold $(M,J,\omega,\theta)$ is called \emph{Vaisman} if the Lee form is parallel with respect to the Levi-Civita connection associated to the metric itself. Vaisman manifolds were intensively studied by Izu Vaisman under the name \emph{generalized Hopf manifolds}. It is known that Hopf manifolds always admits LCK metrics \cite{ornea-verbitsky-JGA}, while they admit Vaisman metrics if and only if they are semisimple, see \cite[Corollary 16.16]{ornea-verbitsky-Book}.
A Vaisman metric is always {\em Gauduchon} in the sense of \cite{gauduchon}, namely, satisfies $\partial\overline\partial\omega^{n-1}=0$, where $n$ is the complex dimension of $M$. Therefore, by \cite[Th\'eor\`eme 1]{gauduchon}, if a Vaisman metric exists, it is unique in its conformal class up to a constant multiplier. Furthermore, any compact Vaisman manifold is LCK with potential \cite[Proposition 4.4]{verbitsky2004vanishing}.
Up to a positive constant, a Vaisman metric is determined by the following identity:
\begin{equation}\label{vaismanomega0det}
\omega = \frac{1}{|\theta|^2}(d^c\theta + \theta \wedge \theta^c),
\end{equation}
where $\theta^c:=J^{-1}\theta J$, and $d^c=J^{-1}dJ$.

Let us explain he role of $\omega_0:=d^c\theta$.
Any compact Vaisman manifold $M$ of LCK rank $1$ is biholomorphically isometric to a Vaisman manifold obtained as $C(S)/\Z$, where $S$ is Sasakian, $C(S) := S \times \R^{>0}$ is its K\"ahler cone, and $\mathbb Z$ is the group generated by $(x,t) \mapsto (\varphi(x), q t)$ for some $q < 1$ and $\varphi$ a Sasakian automorphism of $S$, see \cite{ornea-verbitsky-MRL}. Moreover, the triple $(S, \varphi, q)$ is unique.
We recall that on a Vaisman manifold $M$ the \emph{Lee field} $\theta^\#$, defined by $\theta = g(\theta^\#, \cdot)$, and the \emph{anti-Lee field} $J\theta^\#$ generate a holomorphic totally geodesic foliation $\Sigma$, called the \emph{canonical foliation} \cite{vaisman1979locconf}. We say that a Vaisman manifold $M$ is \emph{quasi-regular} if each leaf of $\Sigma$ is compact, and \emph{irregular} otherwise. Moreover, $M$ is called \emph{regular} if the leaves of $\Sigma$ are orbits of the group $T^2 = (S^1)^2$ freely acting on $M$. We can describe the tangent to $\Sigma$ as the kernel of the positive closed $(1,1)$-form $\omega_0$ so that $(M,\Sigma,\omega_0)$ is a transversally K\"ahler foliation. The projection $\pi \colon M \to \mathcal{X}:=M/\Sigma$ to the leaf space is called the \emph{Lee fibration}. Notice that $\mathcal{X}$ has the structure of a manifold, respectively an orbifold, if the foliation $\Sigma$ is regular, respectively quasi-regular. In such cases there exists a K\"ahler form (respectively K\"ahler orbifold-form) $\omega_{\mathcal{X}}$ on $\mathcal{X}$ such that $\omega_0 = \pi^*\omega_{\mathcal{X}}$ \cite[Proposition 6.4]{verbitsky2004vanishing}.
The structure of regular and quasi-regular Vaisman manifolds will be described in more details in the following sections.

\subsection{Deformations of Vaisman structures}
\label{sec:deformations}

We recall here some deformations of Vaisman structures of LCK rank $1$, which preserve the canonical foliation and the LCK rank.
Deformations of Vaisman structures have been introduced by Tsukada \cite{tsukada1994holomorphic} and recently Madani, Moroianu and Pilca gave an exhaustive presentation in \cite[Section 3]{madani-moroianu-pilca}.
From now on, we will consider $(M, J, \omega, \theta)$ a compact quasi-regular Vaisman manifold of LCK rank $1$. We will always assume that the Vaisman structure is normalized, {\itshape i.e.} the Lee form has norm $1$. Such assumption is not restrictive, since every Vaisman metric is homothetic to a normalized Vaisman metric.

The first kind of deformations we consider are {\em $\Sigma_a$-homotheties}, where $a\in \mathbb R$ is a positive number. They are defined as
\begin{equation*}
\theta_a = a\theta, \quad
\omega_a = a\omega + (a^2-a)\theta \wedge \theta^c.
\end{equation*}
Equivalently, the Vaisman structure $(\omega_a, \theta_a)$ on $M$ can be obtained from the K\"ahler cone by setting the new coordinate $\widetilde{t} = t^a$.
Furthermore, notice that the metric $(\omega_a, \theta_a)$ still has LCK rank $1$, as $\theta_a$ is proportional to $\theta$.

More in general, one can consider \emph{deformations of type I} by substituting, in Equation \eqref{vaismanomega0det}, a new Lee form $\theta' = a\theta + \alpha$, where $\alpha$ is a harmonic $1$-form orthogonal to $\theta$ and $\theta^c$.
In \cite[Lemma 3.3]{madani-moroianu-pilca}, it is proven that the pair $(\omega' := d^c\theta' + \theta' \wedge \theta'^c, \theta')$ is a normalized Vaisman structure on $(M,J)$.

The other family of deformations we consider is given by {\em transverse K\"ahler deformations}, also known as \emph{deformations of type II} \cite{madani-moroianu-pilca}.
Following the same notation as above, denote by $C(S)$ the K\"ahler cone, by $\mathcal X$ the leaf space orbifold of the canonical foliation $\Sigma$ endowed with the K\"ahler orbifold-form $\omega_{\mathcal X}$, and by $\sfrac{1}{4}t^2$ the K\"ahler potential of the initial Vaisman structure. Let $f$ be a basic function with respect to the canonical foliation $\Sigma$ such that $\widetilde{\omega}_{\mathcal{X}} = \omega_{\mathcal{X}} + i\partial\overline{\partial}f > 0$. A transverse K\"ahler deformation is given by replacing $t$ with $\widetilde{t} = e^ft$ and leaving $C(S)$ unchanged.
As we will explain in \Cref{sec:quasi-regular}, the Vaisman structure is determined by the transverse K\"ahler form $\omega_{\mathcal X}$ and the action of $\mathbb Z$ on the K\"ahler cone.
The corresponding K\"ahler form on $C(S)$ is $\widetilde{\Omega} = \sfrac{\sqrt{-1}}{2}\partial\overline{\partial}\widetilde{t}^2 = \sfrac{\sqrt{-1}}{2}\partial\overline{\partial}(e^{2f}t^2)$. The Lee forms are related by $\widetilde{\theta} = \theta -2df$, thus $[\widetilde{\theta}] = [\theta] \in H^1(M,\R)$ so that the associated LCK metric $\widetilde{\omega} := \widetilde{t}^{-2}\widetilde{\Omega}$ has LCK rank $1$, too.
Since the structure on the K\"ahler cone $C(S)$ does not change under the deformation above, the group $\Z$ still acts by biholomorphisms on it. Thus, the metric $\widetilde{\omega}$ is Vaisman.
Notice that $\widetilde{\omega}_{\mathcal{X}}$ defines the same K\"ahler class as $\omega_{\mathcal{X}}$.

\subsection{Approximation of Vaisman structures of generic LCK rank}
\label{approximanion-rank}

Let $(\omega, \theta)$ be a Vaisman structure of LCK rank $k$ on a compact manifold $M$. Following the proof of \cite[Proposition 3.3]{ornea-verbitsky-JGP} by Ornea and Verbitsky, we consider a sequence of cohomology classes $[\alpha_t] \in H^1(M,\R)$ such that $a_t[\theta] + [\alpha_t]$ is rational and approximates $[\theta]$. We take $\alpha_t$ harmonic and pointwise orthogonal to $\theta$ and $J\theta$ so that they define a deformation $(\omega_t, \theta_t)$ of type I.
By construction, $\omega_t$ has LCK rank $1$.
In \cite{ornea-verbitsky-JGP}, it is noted that any metric $\omega_t$ is conformally equivalent to a Vaisman metric, which is in particular Gauduchon.
In fact, the metric $\omega_t$ is already Vaisman itself, as proven in \cite[Lemma 3.3]{madani-moroianu-pilca}.

In conclusion, any Vaisman metric of LCK rank $k$ can be smoothly approximated by Vaisman metrics of LCK rank $1$.

\section{Proof of \texorpdfstring{\Cref{thm:a}}{Theorem A}}
\label{sec:regular}

In this section we prove the approximation theorem in the regular case. Let $M$ be a compact regular Vaisman manifold, with Lee form $\theta$.
As we recalled in \Cref{approximanion-rank}, we assume the LCK rank equal to $1$, without loss of generality.
The structure of $M$ is described in \cite{ornea-verbitsky-MRL}, see also \cite[Theorem 11.6]{ornea-verbitsky-Book}. It is shown that $M$ is biholomorphic to a smooth elliptic fibration $\pi \colon M \to X$ over a projective manifold $X$, which is obtained as the leaf space of the canonical foliation $\Sigma$. More precisely, the manifold $M$ is the quotient of the total space $\widetilde{M} := L^{-1}\setminus\{0\}$ of a negative line bundle $L^{-1}$ over $X$ with the zero section removed. The deck transformation group is $\mathbb Z$ thanks to the assumption of LCK rank being equal to $1$.

The line bundle $L$ is constructed as the pushforward of the complexified weight bundle $L_\C$ over $M$, via $\pi\colon M \to X$.
More precisely, we consider the weight bundle $L_\rho$ of $M$, namely the line bundle associated to the character $\rho$ of the deck transformation group. Consider $L_\rho$ equipped with the flat Weyl connection $\nabla^{\mathcal{W}}$, namely the unique torsion-free connection such that $\nabla^{\mathcal W}g=\theta\otimes g$, see \cite[Theorem 2.2]{vaisman1980locally}. Take the complexified weight bundle $L_\C := L_\rho \otimes_\R \C$ with $\nabla^{\mathcal{W},\C}$ the induced connection. The $(0,1)$-part of $\nabla^{\mathcal{W},\C}$ gives a holomorphic structure on $L_\C$.
Let $L := \pi_*L_\C$ be the corresponding sheaf, in fact line bundle since flat, over $X$.
The space $L^{-1}\setminus\{0\}$ is identified as the K\"ahler cone $C(S)$ over a Sasakian manifold $S$.
Then $M$ is obtained as the quotient of $L^{-1}\setminus\{0\}$ by the action of $\Z \simeq \langle \varphi \times q\rangle$, where $\varphi$ is a finite order Sasakian automorphism of $S$ and $q < 1$.

Fix a trivializing section $\ell$ for $L_{\mathbb C}$.
Consider a Hermitian structure $h$ on $L_\C$ such that $|\ell|_h = 1$. Let $D$ be its Chern connection on $L_\C$.
Since the natural map $L_\C \to \pi^*\pi_*L_\C$ is an isomorphism \cite[Theorem 3.4]{ornea-verbitsky-MathAnn}, $L$ is equipped with a natural Hermitian metric $h_X$, such that the metric $h$ on $L_\C = \pi^*L$ is the pullback of $h_X$.
Therefore, the Chern connection $D$ of $(L_\C,h)$ is lifted from the Chern connection of $(L,h_X)$.
The curvature of $D$ is $-2\sqrt{-1}\omega_0$, where $\omega_0$ is the transverse K\"ahler form, see \cite[Theorem 6.7]{verbitsky2004vanishing}. 
Then, let $-2\sqrt{-1} \omega_X$ be the curvature of $(L,h_X)$, so that $\omega_X$ is a K\"ahler form on $X$ such that $\pi^*\omega_X=\omega_0$.
In particular, $L$ is a positive line bundle over $X$, with $c_1(L)=\sfrac{1}{\pi}[\omega_X]>0$.
We can now describe the Vaisman structure in terms of $h$, $\varphi$ and $q$.
The dual metric $h^{-1}$ on $L^{-1}$ defines the second coordinate $(s,t) \in S \times \R^{>0} = L^{-1}\setminus\{0\}$ by
$t \colon L^{-1}\setminus\{0\} \to \R^{>0}$, $(p,v) \mapsto |v|_{h^{-1}_p}$.
The K\"ahler form on the K\"ahler cone $C(S)$ is given by $\Omega = \sfrac{\sqrt{-1}}{2}\partial\overline{\partial}(t^2)$.
The K\"ahler metric on the projective manifold $X$ is given by $\omega_X = -\sfrac{\sqrt{-1}}{2} \partial\overline{\partial}\log h$.
The Vaisman structure can be read from these data as 
\begin{equation}\label{eq:vaisman}
\omega=t^{-2}\Omega, \quad \theta = -d\log (t^2), \quad \Z \simeq \langle \varphi \times q \rangle \curvearrowright L^{-1}\setminus\{0\}.
\end{equation}

Now, by the Kodaira Embedding Theorem, there exists $k\in\mathbb N$ sufficiently large such that the Kodaira map $\psi_k \colon X \hookrightarrow \CP^{N_k}$ is an embedding, where $N_k + 1 = \dim_{\C}H^0(X,L^k)$. Recall that  $\psi_k^*\mathcal{O}_{\CP^{N_k}}(-1) = L^{-k}$. Thus there is a well defined holomorphic embedding $\widetilde{\psi}_k \colon L^{-k}\setminus\{0\} \hookrightarrow \mathcal{O}_{\CP^{N_k}}(-1)\setminus\{0\} \simeq \C^{N_k +1}\setminus\{0\} $ making the following diagram commute.
\begin{center}\begin{tikzcd}
L^{-k}\setminus\{0\} \arrow[hookrightarrow,r,"\widetilde{\psi}_k"]\arrow[d] &\mathcal{O}_{\CP^{N_k}}(-1)\setminus\{0\}\arrow[d]\\
X\arrow[hookrightarrow,r,"\psi_k"] & \CP^{N_k} \end{tikzcd} \end{center}

The holomorphic map $p_k \colon L^{-1}\setminus\{0\} \to L^{-k}\setminus\{0\}$ mapping $v$ to $v^{\otimes k}$ is a $k$-fold covering. Thus, we can define a holomorphic immersion \begin{align*}\phi_k := \widetilde{\psi}_k \circ p_k \colon L^{-1}\setminus\{0\} \to \C^{N_k + 1} \setminus\{0\}.\end{align*}

Consider now the action of $\mathbb Z$ on $L^{-1}\setminus\{0\}$. Note that the generator $\gamma := \varphi \times q$ of $\mathbb Z$ commutes with the $\C^*$-action generated by the Reeb flow of $S$ and the action of $\R^{>0}$, hence $\gamma$ yields an automorphism of the line bundle $L$.
Moreover, recall that $\gamma$ acts as a contraction.
Therefore, $\gamma$ can be extended to a contraction $\Gamma$ of the line bundle $\mathcal{O}_{\CP^{N_k}}(-1)\setminus\{0\}$ such that
\begin{equation*} \phi_k \circ \gamma = \Gamma \circ \phi_k.\end{equation*}
Summarizing, we have that $\phi_k$ descends to a holomorphic immersion
$$ \varphi_k \colon M = ({L^{-1}\setminus\{0\}}) / {\langle \gamma \rangle} \rightarrow ({\mathcal{O}_{\CP^{N_k}}(-1)\setminus\{0\}}) / {\langle \Gamma \rangle} =: \mathcal{H}_{N_k + 1} $$
of $M$ into a Hopf manifold $\mathcal{H}_{N_k + 1}$ of dimension $N_k + 1$, making the following diagram commutative.
\begin{center}\begin{tikzcd} M \arrow[rightarrow,r,"\varphi_k"]\arrow[d,"\pi"] &\mathcal{H}_{N_k+1}\arrow[d,"\pi_k"]\\
X\arrow[hookrightarrow,r,"\psi_k"] & \CP^{N_k} \end{tikzcd}\end{center}
The Hopf manifold $\mathcal{H}_{N_k + 1}$ is Vaisman by construction and regular, since $\CP^{N_k}$ is a projective manifold.
In particular, it is also semisimple, see \cite[Corollary 16.16]{ornea-verbitsky-Book} (for complex surfaces, see \cite{gauduchon-ornea}).

Now, as in \cite[Proposition 5.1]{loi-placini-zedda}, we notice the following.
Fix an orthonormal basis $\{s_0, \dots, s_{N_k}\}$ of $H^0(X,L^k)$ with respect to the standard $L^2$-product induced by $h_X$ and $\omega_X$.
The immersion $\phi_k \colon L^{-1} \setminus \{0\} \to H^0(X,L^k)^*$ is then given by
$$ \phi_k (v) = \left( s \mapsto \frac{ v^k(s) }{\|s\|_{L^2}} \right) .$$
We can then compute the pullback via $\phi_k$ of the Fubini-Study metric $h^{-1}_{FS}$ on $\mathcal{O}_{\CP^{N_k}}(-1)$, for $v\in L^{-1}\setminus\{0\}$, as
\begin{align*}\phi_k^*(h_{FS}^{-1})_x(v,v) = h^{-k}(v,v)\sum_{j=0}^{N_k}h(s_j(x),s_j(x)) 
= \B_k(x)h^{-k}(v,v),\end{align*} where $\B_k(x) = \sum\limits_{j=0}^{N_k}h(s_j(x),s_j(x))$ is the Bergman kernel of $(X,L^k)$. 
Note that $\phi_k^*(h_{FS}^{-1})$ is not a Hermitian metric on the line bundle $L^{-1}$ since it does not scale correctly under the $\C^*$-action. 
Nevertheless, it defines a change of the coordinate on $\R^{>0}$ of the K\"ahler cone $\widetilde{M}$, corresponding to the composition of the $\Sigma_k$-homothetic transformation $t \mapsto t^k$, with the transverse K\"ahler deformation $t^k \mapsto \B_kt^k$.

We have proved the following immersion result:
\begin{proposition}\label{immersionregularVaisman}
Let $M$ be a compact regular Vaisman manifold of LCK rank $1$.
For every sufficiently large positive integer $k$, there exists a holomorphic immersion $\varphi_k \colon M \to \mathcal{H}_{N_k+1}$ into a semisimple regular Hopf manifold of dimension $N_k+1$.
Moreover, if we denote by $\phi_k\colon  \widetilde M \to \mathbb C^{N_k+1}\setminus\{0\}$ the immersion lifting $\varphi_k$ to the K\"ahler coverings and by $t$ and $\tau$ the coordinates on the cone direction $\mathbb R^{>0}$ of $M$ and $\mathcal{H}_{N_k+1}$ respectively, then $\phi_k^*(\tau) = \B_kt^k$.
\end{proposition}

We are now able to conclude the proof of \Cref{thm:a} by using the sequence of holomorphic immersions $\varphi_k \colon M \to \mathcal{H}_{\N^k+1}$ constructed above.
Notice that $\tau$ is the coordinate induced by the flat metric on $\C^{N_k + 1}\setminus\{0\}$ or, equivalently, by the Hermitian metric $h_{FS}$ on $\mathcal{O}_{\CP^{N_k}}(-1)\setminus\{0\}$ whose curvature is $-\omega_{FS}$.

The holomorphic immersions $\varphi_k \colon M \to \mathcal{H}_{N_k+1}$ induce Vaisman structures on $M$ determined by $\phi_k^*(\tau),\varphi$ and $q$ as in Equation \eqref{eq:vaisman}.
Now the $\Sigma_{\sfrac{1}{k}}$-homothety of these structures are transverse K\"ahler deformations, determined by the Bergman kernel $\B_k$, of the original Vaisman structure. By the Tian approximation, the first coefficient of the asymptotic expansion of the Bergman kernel $\B_k$ smoothly converges to $1$ when $k$ goes to infinity. Summarizing, the $\Sigma_{\sfrac{1}{k}}$-homotheties of the induced structures converge smoothly to $(\omega,\theta)$.

\section{Proof of \texorpdfstring{\Cref{thm:b}}{Theorem B}}
\label{sec:quasi-regular}

In this section we prove the approximation theorem in the quasi-regular case. Let $M$ be a compact quasi-regular Vaisman manifold, with Lee form $\theta$.
As explained in \Cref{approximanion-rank}, there is no loss of generality in assuming the LCK rank equal to $1$. We want to construct a holomorphic embedding of $M$ into a semisimple quasi-regular Hopf manifold.
Although the proof is analogous to the regular case, we repeat the main steps for the sake of clarity.

We recall the structure of $M$, see \cite[Theorem 11.9]{ornea-verbitsky-Book}.
We can describe $M$ as an elliptic fibration $\pi \colon M \to \mathcal{X}$ over a projective orbifold $\mathcal X$, obtained as the leaf space of the canonical foliation $\Sigma$. 
As we are going to explain, the orbifold $\mathcal X$ is endowed with a positive line orbibundle $L$, with curvature $\omega_{\mathcal X}>0$.
The K\"ahler covering $\widetilde M$ of $M$ is identified with $L^{-1}\setminus\{0\}$, the total space of the negative line orbibundle $L^{-1}$ without the zero section. It is also identified with the K\"ahler cone $C(S)$ of a compact quasi-regular Sasakian manifold $S$.
Therefore the manifold $M$ is the quotient of $L^{-1}\setminus\{0\}$ by the action of $\Z \simeq \langle \varphi \times q\rangle$, where $\varphi$ is a finite order Sasakian automorphism of $S$ and $q < 1$. The line orbibundle $L$ is given by the pushforward of the complexified weight bundle $(L_\C, h)$ over $M$, where $h$ is a Hermitian metric on $L_\C$ so that the trivializing section $\ell$ has norm $|\ell|_h = 1$.
Under this identification, the induced Hermitian metric $h_{\mathcal{X}}$ on $L$ has Chern curvature $-2\sqrt{-1}\omega_{\mathcal{X}}$. This shows that $L$ is an ample line orbibundle.
In analogy with the regular case, the Vaisman structure is described in terms of $h$, $\varphi$ and $q$, see Equation \eqref{eq:vaisman}, keeping into account that the objects are to be regarded in the orbifold sense.

Since $L$ is an ample line orbibundle, for any fixed $m\in\mathbb N$ and for $k\in\mathbb N$ sufficiently large, Ross and Thomas \cite{ross-thomas-2} defined an orbifold embedding $\psi_k \colon \mathcal{X} \hookrightarrow \CP^{N_k}(\textbf{w})$ of $\mathcal{X}$ into the weighted projective space $\CP^{N_k}(\textbf{w})$, for a suitable weight $\textbf{w}$.
Furthermore, it holds \begin{align*}\psi_k^*\mathcal{O}_{\CP^{N_k}(\textbf{w})}(-1) = L^{-1}.\end{align*}
We have the following commutative square, where the horizontal arrows are embeddings:
\begin{center}
\begin{tikzcd}L^{-1}\setminus\{0\}\arrow[hookrightarrow,r,"\widetilde{\psi}_k"]\arrow[d] &\mathcal{O}_{\CP^{N_k}(\textbf{w})}(-1)\setminus\{0\}\arrow[d]\\
\mathcal{X}\arrow[hookrightarrow,r,"\psi_k"] & \CP^{N_k}(\textbf{w}) . \end{tikzcd}
\end{center}

We proceed as for the regular case: the action of $\gamma=\varphi \times q$ on $L^{-1}\setminus\{0\}$ commutes with the $\C^*$-action generated by the Reeb flow and the action of $\R^{>0}$. Hence $\gamma$ yields an automorphism of the line orbibundle $L$.
Moreover, $\gamma$ can be extended to a contraction $\Gamma$ of the line orbibundle $\mathcal{O}_{\CP^{N_k}(\textbf{w})}(-1)\setminus\{0\}$ such that 
\begin{equation*} \widetilde{\psi}_k \circ \gamma = \Gamma \circ \widetilde{\psi}_k.\end{equation*}
Again, $\widetilde{\psi}_k$ descends to a holomorphic embedding
$$\varphi_k \colon M = ({L^{-1}\setminus\{0\}}) / {\langle \gamma \rangle} \hookrightarrow ({\mathcal{O}_{\CP^{N_k} (\textbf{w})}(-1)\setminus\{0\}})/{\langle \Gamma \rangle} =: \mathcal{H}_{N_k + 1} $$
of $M$ into a Hopf manifold $\mathcal{H}_{N_k + 1}$ of dimension $N_k + 1$, making the following diagram commutative: 
\begin{center}\begin{tikzcd} M \arrow[hookrightarrow,r,"\varphi_k"]\arrow[d,"\pi"] &\mathcal{H}_{N_k+1}\arrow[d,"\pi_k"]\\
\mathcal{X}\arrow[hookrightarrow,r,"\psi_k"] & \CP^{N_k}(\textbf{w}) . \end{tikzcd}\end{center} 
The Hopf manifold $\mathcal{H}_{N_k + 1}$ is Vaisman by construction and quasi-regular, since $\CP^{N_k}(\textbf{w})$ is a projective orbifold.
In particular, it is also semisimple, see \cite[Corollary 16.16]{ornea-verbitsky-Book}.

We have proved the following embedding result:
\begin{proposition}\label{embeddingquasiregularVaisman} Let $(M,\omega,\theta)$ be a compact quasi-regular Vaisman manifold of LCK rank $1$.
For every sufficiently large positive integer $k$, there exists a holomorphic embedding $\varphi_k \colon M \to \mathcal{H}_{N_k+1}$ into a semisimple quasi-regular Hopf manifold.
\end{proposition}

We are now able to prove the approximation \Cref{thm:b}.
For a fixed $m\in\mathbb N$, we have a sequence of holomorphic embeddings $\varphi_k \colon M \to \mathcal{H}_{N_k+1}$ such that $\psi_k^*(\mathcal{O}_{\CP^{N_k}(\textbf{w})}(-1)) = L^{-1}$.
Moreover, the Vaisman structure of $\mathcal H_{N_k+1}$ is determined by the action of $\Gamma$ and the Hermitian metric $h'_{FS}$ (following the notation in \cite{loi-placini}) on $\mathcal{O}_{\CP^{N_k}(\textbf{w})}(-1)$ whose curvature is the Fubini-Study form $\omega_{FS}$ of the weighted projective space $\CP^{N_k}(\textbf{w})$.
By \cite[Proposition 17]{loi-placini}, the sequence $\psi_k^*h'_{FS}$ converges in $\mathcal C^m$-norm to $h$.
Therefore, the Vaisman structures on $M$ determined by $\psi_k^*h'_{FS}$, $\varphi$ and $q$ converge in $\mathcal C^m$-norm to $(\omega,\theta)$.

\bibliographystyle{alpha}
\bibliography{biblio}

\begin{thebibliography}{MMP21}

\bibitem[AZ19]{angella2019isometric}
Daniele Angella and Michela Zedda.
\newblock Isometric immersions of locally conformally {K}{\"a}hler manifolds.
\newblock {\em Annals of Global Analysis and Geometry}, 56(1):37--55, 2019.

\bibitem[Bai57]{baily}
W.~L. Baily.
\newblock On the imbedding of {$V$}-manifolds in projective space.
\newblock {\em Amer. J. Math.}, 79:403--430, 1957.

\bibitem[BG08]{boyer-galicki-Book}
Charles~P. Boyer and Krzysztof Galicki.
\newblock {\em Sasakian geometry}.
\newblock Oxford Mathematical Monographs. Oxford University Press, Oxford,
  2008.

\bibitem[Gau77]{gauduchon}
Paul Gauduchon.
\newblock Le th\'eor\`eme de l'excentricit\'e{} nulle.
\newblock {\em C. R. Acad. Sci. Paris S\'er. A-B}, 285(5):A387--A390, 1977.

\bibitem[GO98]{gauduchon-ornea}
P.~Gauduchon and L.~Ornea.
\newblock Locally conformally {K}\"ahler metrics on {H}opf surfaces.
\newblock {\em Ann. Inst. Fourier (Grenoble)}, 48(4):1107--1127, 1998.

\bibitem[LP22]{loi-placini}
Andrea Loi and Giovanni Placini.
\newblock Any sasakian structure is approximated by embeddings into spheres,
  2022.
\newblock to appear in {\itshape Forum Math.} \texttt{arXiv:2210.00790}.

\bibitem[LPZ24]{loi-placini-zedda}
A.~Loi, G.~Placini, and M.~Zedda.
\newblock Immersions into {S}asakian space forms.
\newblock {\em Math. Z.}, 307(3):Paper No. 60, 2024.

\bibitem[MM07]{ma-marinescu-Book}
Xiaonan Ma and George Marinescu.
\newblock {\em Holomorphic {M}orse inequalities and {B}ergman kernels}, volume
  254 of {\em Progress in Mathematics}.
\newblock Birkh\"auser Verlag, Basel, 2007.

\bibitem[MMP21]{madani-moroianu-pilca}
Farid Madani, Andrei Moroianu, and Mihaela Pilca.
\newblock Lc{K} structures with holomorphic {L}ee vector field on
  {V}aisman-type manifolds.
\newblock {\em Geom. Dedicata}, 213:251--266, 2021.

\bibitem[OV03]{ornea-verbitsky-MRL}
Liviu Ornea and Misha Verbitsky.
\newblock Structure theorem for compact {V}aisman manifolds.
\newblock {\em Math. Res. Lett.}, 10(5-6):799--805, 2003.

\bibitem[OV05]{ornea-verbitsky-MathAnn}
Liviu Ornea and Misha Verbitsky.
\newblock An immersion theorem for {V}aisman manifolds.
\newblock {\em Math. Ann.}, 332(1):121--143, 2005.

\bibitem[OV10]{ornea2010locally}
Liviu Ornea and Misha Verbitsky.
\newblock Locally conformal {K}{\"a}hler manifolds with potential.
\newblock {\em Mathematische Annalen}, 348(1):25--33, 2010.

\bibitem[OV16]{ornea-verbitsky-JGP}
Liviu Ornea and Misha Verbitsky.
\newblock L{CK} rank of locally conformally {K}\"{a}hler manifolds with
  potential.
\newblock {\em J. Geom. Phys.}, 107:92--98, 2016.

\bibitem[OV23]{ornea-verbitsky-JGA}
Liviu Ornea and Misha Verbitsky.
\newblock Non-linear {H}opf manifolds are locally conformally {K}\"ahler.
\newblock {\em J. Geom. Anal.}, 33(7):Paper No. 201, 10, 2023.

\bibitem[OV24]{ornea-verbitsky-Book}
Liviu Ornea and Misha Verbitsky.
\newblock {\em Principles of locally conformally {K{\"a}hler} geometry}.
\newblock Cham: Birkh{\"a}user, 2024.

\bibitem[RT11a]{ross-thomas-1}
Julius Ross and Richard Thomas.
\newblock Weighted {B}ergman kernels on orbifolds.
\newblock {\em J. Differential Geom.}, 88(1):87--107, 2011.

\bibitem[RT11b]{ross-thomas-2}
Julius Ross and Richard Thomas.
\newblock Weighted projective embeddings, stability of orbifolds, and constant
  scalar curvature {K}\"{a}hler metrics.
\newblock {\em J. Differential Geom.}, 88(1):109--159, 2011.

\bibitem[Rua98]{ruan1998canonical}
Wei-Dong Ruan.
\newblock Canonical coordinates and {B}ergmann metrics.
\newblock {\em Communications in Analysis and Geometry}, 6(3):589--631, 1998.

\bibitem[Tia90]{tian1990set}
Gang Tian.
\newblock On a set of polarized {K}{\"a}hler metrics on algebraic manifolds.
\newblock {\em Journal of Differential Geometry}, 32(1):99--130, 1990.

\bibitem[Tsu94]{tsukada1994holomorphic}
Kazumi Tsukada.
\newblock Holomorphic forms and holomorphic vector fields on compact
  generalized hopf manifolds.
\newblock {\em Compositio Mathematica}, 93(1):1--22, 1994.

\bibitem[Vai79]{vaisman1979locconf}
Izu Vaisman.
\newblock Locally {C}onformal {K}{\"a}hler {M}anifolds with {P}arallel {L}ee
  {F}orm.
\newblock {\em Rend. Mat. Roma}, 12:263–--284, 1979.

\bibitem[Vai80]{vaisman1980locally}
Izu Vaisman.
\newblock On locally and globally conformal {K}{\"a}hler manifolds.
\newblock {\em Transactions of the American Mathematical Society},
  262(2):533--542, 1980.

\bibitem[Ver04]{verbitsky2004vanishing}
Mikhail~Sergeevich Verbitsky.
\newblock Vanishing theorems for locally conformal hyperk{\"a}hler manifolds.
\newblock {\em Trudy Matematicheskogo Instituta imeni VA Steklova}, 246:64--91,
  2004.

\bibitem[Zel98]{zelditch1998szego}
Steve Zelditch.
\newblock Szeg{\"o} kernels and a {T}heorem of {T}ian.
\newblock {\em International Mathematics Research Notices}, 1998(6):317--331,
  1998.

\end{thebibliography}

\end{document}